\documentclass[authoryear]{elsarticle}
\usepackage[latin9]{inputenc}
\usepackage[T1]{fontenc}
\usepackage{epsfig}
\usepackage{amsmath}
\usepackage{listings}
\usepackage[a4paper]{geometry}
\begin{document}
   \title{A simple MATLAB program to compute differentiation matrices for arbitrary meshes via Lagrangian interpolation}
    \author{Miguel P\'erez-Saborid\corref{cor1}}
   \address{Departamento de Ingenier\'ia Aeroespacial y Mec\'anica de Fluidos, Escuela T\'ecnica Superior de Ingenier\'ia, Universidad de Sevilla, Av. de los Descubrimientos s/n, 41092 Sevilla, Spain}

 \ead{psaborid@gmail.com}
 \cortext[cor1]{The author will appreciate any comments or suggestions to improve later versions of this document.  He
would also like to thank Prof. Bosco Garc\'ia-Archilla for his careful and insightful reading of the manuscript.}

\begin{abstract}
A MATLAB program for computing differentiation matrices for arbitrary one-dimensional meshes is presented in this manuscript. The differentiation matrices for a mesh of $N$ arbitrarily spaced points are formed from
those obtained using Lagrangian interpolation on stencils of a fixed but arbitrary number $M\le N$ of contiguous mesh points. For the particular case $M=N$ and meshes with  Chebyshev or Legendre distributions of points, the program yields  the well known spectral differentiation matrices. For $M<N$ and $M$ odd, the differentiation matrices coincide, for the special case of an evenly spaced mesh,  with those obtained by central finite differences. 
\end{abstract}

\maketitle

\section{Introduction}
The Lagrangian interpolation approach to differentiation matrices has the advantage over that based on Taylor series expansions that it does not yield poorly condtioned Vandermonde systems of equations when the number of stencil points is large [Leveque (2007), Fornberg (1996)]. In this brief document, a basic MATLAB program for the computation of differentiation matrices for arbritrary meshes using Lagrangian interpolation is presented. As will be seen, the approach allows a systematic use of the differentiation matrices for stencils of $M$ arbitrarily spaced points to construct those for arbitrarily spaced meshes of $N\ge M$ points. The stencil differentiation matrices are obtained in Section $2$ by means of a known recursion relation for the derivatives of the Lagrangian interpolants [Fornberg (1996), Welfert (1997)] together with a reordering of the quotients of products of distances between stencil points appearing in that relation. This reordering is introduced in order to alleviate round-off error problems associated with the computations of those terms for stencils whose points are separated by widely  varying distances. For the particular case $M=N$ - i.e., when the mesh and the stencil coincide - the differentiation matrices reproduce the well known spectral ones for meshes with Chebyshev or Legendre distributions of points [Trefethen (2001)]. However, sparser differentiation matrices obtained for $M < N$  (usually $M<<N$) are desired in many applications, specially those involving higher dimensions. Thus, it is shown in Section $3$ how to use in these cases the stencil differentiation matrices to construct those for an arbitrarily spaced, larger mesh. In the particular case of an evenly spaced mesh, the method proposed in this manuscript yields the differentiation matrices corresponding to central finite differences of second ($M=3$) and higher orders ($M>3$ and odd) [Fornberg (1996)]. Only one-dimensional meshes have been considered in this manuscript, since differentiation matrices for higher dimensions can be obtained from one-dimensional ones through the well known matrix Kronecker product [Trefethen (2001)].
 
\section{The stencil}
The Lagrangian interpolation formula for a function $f(x)$ in an interval containing a stencil of $M$ points \footnote{The points and differentiation matrices corresponding to a stencil will be denoted by an overbar to distinguish them from those corresponding to the general mesh introduced in the next section.} $\bar{x}_{1}$, $\bar{x}_{2}$,...$\bar{x}_{M}$ is  [Henrici (1982)]
\begin{equation}
f(x)=\sum_{m=1}^{M} f_m L_m(x), \label{fint}
\end {equation}
where $\bar{x}_{1} \le x \le \bar{x}_{M}$, $f_m=f(\bar{x}_{m})$ and $L_m(x)$ is the Lagrangian interpolant associated to the stencil point $\bar{x}_{m}$,
\begin{equation}
L_m(x)=\frac{(x-\bar{x}_{1})...(x-\bar{x}_{m-1})(x-\bar{x}_{m+1})....(x-\bar{x}_{M})}{(\bar{x}_{m}-\bar{x}_{1})...(\bar{x}_{m}-\bar{x}_{m-1)})(\bar{x}_{m}-\bar{x}_{m+1)})....(\bar{x}_{m}-\bar{x}_{M})}. \label{Lint}
\end {equation}

 Despite their complicated appearance, the Lagrangian interpolants satisfy beautiful relations that can often be demonstrated by simple arguments. For instance, by making $f(x)=x^r$  in (\ref{fint}) we obtain
\begin{equation}
 \bar{x}_1^r L_1(x)+\bar{x}_2^r L_2(x)+...\bar{x}_M^r L_M(x)=x^r \quad {\rm for}\quad 0 \le r \le M-1,   \label{relation}
\end {equation}
which follows from the fact that  Lagrange's interpolation formula (\ref{fint}) must reproduce exactly any polynomial of degree less or equal than  $M-1$. In particular,  relation (\ref{relation}) yields for $r=0$ that $\sum_{m=1}^M L_m(x) =1$, or
\begin{equation}
\sum_{m=1}^M L'_m(x) =0,    \label{relation1}
\end {equation}
a result which will be very useful in the sequel.

It is often convenient to express the Lagrangian interpolants in the so-called \emph{barycentric form} [Rutishauser (1976), Henrici(1982), Barut and Trefethen (2004)]
\begin{equation}
L_m(x)=\frac{a(x)}{(x-\bar{x}_m) a_m}\quad (m=1,...M), \label{Lbary}
\end {equation}
where the polinomial $a(x)$ and the constant $a_m$ are defined by
\begin{equation}
a(x)\equiv \prod_{k=1}^{M} (x-\bar{x}_{k})\quad {\rm and}\quad a_m\equiv a'(\bar{x}_m) =\prod_{k\ne m} (\bar{x}_m-\bar{x}_{k}). \label{aam}
\end{equation}
The barycentric form is especially useful for finding the derivatives of $L_m(x)$ at any point $\bar{x}_i$ of the stencil. In fact,  by writing (\ref{Lbary}) as
\begin{equation}
a(x)=a_m (x-\bar{x}_m
) L_m (x) \label{Lbary2}
\end{equation}
and taking its derivative,
\begin{equation}
a'(x)=a_m (x-\bar{x}_m) L'_m (x)+a_m L_m (x), \label{Lbaryd}
\end{equation}
one immediately obtains  for any stencil point $\bar{x}_i$ different from $\bar{x}_m$  that  
\begin{equation}
L'_m (\bar{x}_i)=\frac{a_i}{a_m (\bar{x}_i-\bar{x}_m)} \quad ( m \ne i) , \label{Lpjnem}
\end{equation}
where we have taken into account the definition $a'(\bar{x}_i)=a_i$ and that $L_m (\bar{x}_i)=0$ if $m \ne i$. In order to find $L'_i(\bar{x}_i)$, it is more convenient, however, to use (\ref{relation1}) particularized at $\bar{x}_i$ together with (\ref{Lpjnem}), which yield
\begin{equation}
L'_i (\bar{x}_i)=-\sum_{m=1,m\ne i}^M L'_m (\bar{x}_i) =-\sum_{m=1,m\ne i}\frac{a_i}{a_m ( \bar{x}_i-\bar{x}_m)}. \label{Lpjeqm}
\end{equation}

Higher order derivatives of the Lagrangian interpolants can be efficiently computed using a recursion relation [Fornberg (1996), Welfert (1997)] which follows directly from equation (\ref{Lbary2}) after taking its $s$-th derivative, 
\begin{equation}
a^{(s)}(x)=a_m (x-\bar{x}_m) L^{(s)}_m (x) +s a_m  L^{(s-1)}_m (x).  \label{Leibniz}
\end{equation}
Note that when  (\ref{Leibniz}) is particularized at $x=\bar{x}_m$ provides $a^{(s)}(\bar{x}_m)= s a_m  L^{(s-1)}_m (\bar{x}_m)$ for any $m=1,...M$. Thus, at $\bar{x}_i\ne \bar{x}_m$   (\ref{Leibniz}) yields
\begin{equation}
 L^{(s)}_m (\bar{x}_i) = \frac{s}{(\bar{x}_i-\bar{x}_m)} \left[\frac{a_i}{a_m} L^{(s-1)}_i (\bar{x}_i)-L^{(s-1)}_m (\bar{x}_i)\right] \quad (i \ne m).  \label{Lpsjnem}
\end{equation}
The value of  $L^{(s)}_i (\bar{x}_i)$ can again be obtained from (\ref{relation1}) as
\begin{equation}
L^{(s)}_i (\bar{x}_i)=-\sum_{m=1,m\ne i}^M L^{(s)}_m (\bar{x}_i). \label{Lpsjeqm}
\end{equation}
Observe that the case $s=1$, which corresponds to  (\ref{Lpjnem})-(\ref{Lpjeqm}), is included in (\ref{Lpsjnem})-(\ref{Lpsjeqm}) if we define $L^{(0)}_m (\bar{x}_i)\equiv \delta_{im}$, where $\delta_{im}$ is the Kronecker's delta symbol.

Expressions (\ref{Lpjnem})-(\ref{Lpjeqm}), together with relations  (\ref{Lpsjnem})-(\ref{Lpsjeqm}), permit the numerical computation of the $s$-th derivative of any Lagrangian interpolant at any point of the stencil. These derivatives form the so-called s-th differentiation matrix of the stencil, denoted by $\bar{D}_s$, which is defined as that whose elements are
\begin{equation}
\bar{D}_s (i,m)=L^{(s)}_m (\bar{x}_i)\quad (i,m=1,...M). \label{Ds}
\end{equation}
Observe that, according to (\ref{fint})  and  (\ref{Ds}), the value of the $s$-th derivative of $f(x)$ at the stencil point $\bar{x}_i$ can be approximated  using  $\bar{D}_s$ as 
\begin{equation}
f^{(s)}(\bar{x}_i)=\sum_{k=1}^{M} \bar{D}_s (i,m)  f(\bar{x}_m) . \label{fpxi}
\end{equation}
The recursion relation  for the differentiation matrices follows directly from  (\ref{Lpsjnem}) as 
\begin{equation}
 \bar{D}_s (i,m) = \frac{s}{(\bar{x}_i-\bar{x}_m)} \left[\frac{a_i}{a_m} \bar{D}_{s-1} (i,i)-\bar{D}_{s-1} (i,m)\right] \quad (i \ne m),  \label{Dbsjnem}
\end{equation}
while  the determination of their diagonal elements follows from (\ref{Lpsjeqm}) as
\begin{equation}
\bar{D}_s (i,i)=-\sum_{m=1,m\ne i}^M \bar{D}_s (i,m). \label{Dbsjeqm}
\end{equation}

   Listed below is the MATLAB function \emph{\bf Dbs.m} that implements equations  (\ref{Dbsjnem})-(\ref{Dbsjeqm}) in order to obtain the elements of the $i$-th row - i.e., those involved in the derivatives of $f(x)$ at the stencil point $\bar{x}_i$
 - of the $\bar{D}_s$ matrices ($s=1,2,3$). The output of the function is given by the row vectors defined as $\bar{d}_s(1:M)\equiv \bar{D}_s(i,1:M)$ ($s=1,2,3$)\footnote{The author has attempted to write a basic code, so that, although it might not be optimal concerning efficiency, it is understandable and will allow for the translation to other coding languages.}.  Higher order differentiation matrices can be trivially included in the program. In addition,  it should be possible to easily modify the code to account for appropriate weight functions in (\ref{Lint}) and  (\ref{Dbsjnem})-(\ref{Dbsjeqm}) if differentiation matrices obtained by methods other than polynomial interpolation - such as Hermite or Laguerre collocation - are desired [Welfert(1997), Weideman and Reddy (1999)].  It is worth pointing out that, due to machine round-off errors, special care has to be taken for a stable computation of the quotients $a_i/a_m$ appearing in the recursion relation, particularly if the stencil is large and has points separated by widely different distances. In order to try to alleviate this problem, the factors intervening in $a_i$ and $a_m$  are first stored in vectors $F_i(1:M)$ and $F_m(1:M)$ whose components are defined as
\begin{equation}
\begin{split}
F_i (k)=\bar{x}_i-\bar{x}_k \quad (k\ne i) \quad{\rm and} \quad F_i (k=i)=1 \quad\\
 F_m (k)=\bar{x}_m-\bar{x}_k \quad (k\ne  m) \quad {\rm and} \quad  F_m (k=m)=1. \label{Fs}
\end{split}
\end{equation}
Then, vectors $F_i$ and $F_m$ are sorted in ascending order of the absolute values os their components and the quotient $a_i/a_m$ computed as
\begin{equation}
\frac{a_i}{a_m}=(-1)^{n_i-n_m} \frac{|F_i(k_1)|}{|F_m(l_1)|} \frac{|F_i(k_2)|}{|F_m(l_2)|} ...\frac{|F_i(k_M)|}{|F_m(l_M)|}, \label{aiam}
\end{equation}
where $[F_i(k_1), F_i(k_2),...F_i(k_M)]$ and $[F_m(l_1), F_m(l_2),...F_m(l_M)]$ represent the reordered components of vectors $F_i$ and $F_m$ and $n_i$ and $n_m$ are their respective number of negative components.  As shown in the function listing, both the reordering and the term by term product of quotients in (\ref{aiam}) can be carried out by by the Matlab functions \emph{{\bf sort}} and \emph{{\bf prod}}, while $n_i$ and $n_m$ can be found by the function \emph{{\bf find}}. Finally, for a fixed $i$, the values of $a_i/a_m$ for all values of $m\ne i$ are stored in vector \emph{{\bf aiam}}(1:M) to be used in the computation of the differentiation matrices. The letter \emph{b} appended in the function to the stencil points and the derivative matrices stands for the overbar which characterizes these quantities in the text.   

\vspace{5 mm}

\lstset{language=Matlab, breaklines=true, basicstyle=\footnotesize}
\lstset{numbers=left, numberstyle=\tiny, stepnumber=1, numbersep=-2pt}
\begin{lstlisting}[frame=single]

   function[db_1 db_2 db_3]=Dbs(i,xb,M)
   %Compute factors a_i/a_m and store values for each m not equal to i 
   aiam=zeros(1,M);
   F_i=xb(i)-xb; F_i(i)=1; n_i=length(find(F_i<0));
   for m=1:M,
    if abs(m-i)>0,
    F_m=xb(m)-xb; F_m(m)=1; n_m=length(find(F_m <0));     
    aiam(m)=(-1)^(n_i-n_m)*prod(sort(abs(F_i))./sort(abs(F_m)));
    end
   end
  %
  % Differentiation matrices:
  Db_1=zeros(M,M); Db_2=zeros(M,M); Db_3=zeros(M,M);
  % s=1:       
  for m=1:M
    if abs(m-i)>0,
    Db_1(i,m)=aiam(m)/(xb(i)-xb(m));
    end
  end     
  Db_1(i,i)=-sum(Db_1(i,:));
  % s=2
  for m=1:M
    if abs(m-i)>0,
     Db_2(i,m)=2/(xb(i)-xb(m))*(aiam(m)*Db_1(i,i)-Db_1(i,m));    
    end
  end     
  Db_2(i,i)=-sum(Db_2(i,:));
  % s=3
  for m=1:M
    if abs(m-i)>0,
     Db_3(i,m)=3/(xb(i)-xb(m))*(aiam(m)*Db_2(i,i)-Db_2(i,m));    
    end
  end     
  Db_3(i,i)=-sum(Db_3(i,:));
  %
  % Output
  db_1=Db_1(i,:); db_2=Db_2(i,:); db_3=Db_3(i,:);
\end{lstlisting}

\pagebreak

\section{Differentiation matrices for an arbitrary mesh}
Consider a function $f(x)$ defined on an arbitrary mesh of $N$ points, $x_1,...x_i,...x_N$, by its values  $f(x_i)\equiv f_i$. A differentiation matrix for the mesh, $D_s$, is one which permits to approximate the values  of the $s$-th derivative of $f$ at the mesh points, $f^{(s)}(x_i)\equiv f^{(s)}_i$, in term of those f as
\begin{equation}
f^{(s)}_i=\sum_{j=1}^N D_s (i, j) f_j\quad (i=1,... N). \label{fpi}
\end{equation}  
A convenient way to find the matrix elements $D_s(i, j=1,...M)$ corresponding  to a  given point $x_i$ is to define a stencil of $M \le N$ contiguous mesh points $(x_{j_1}, x_{j_2}, ...x_{j_M})$  surrounding  $x_i$, to store their coordinates in the vector ${\bf \bar{x}}=\left(\bar{x}_{1}, \bar{x}_{2},...\bar{x}_{M}\right)$ and to use the stencil differentiation matrices obtained in the previous section with the MATLAB function  \emph{{\bf Dbs.m}}.  Thus,  the mesh point $x_i$ coincides with a certain stencil point $\bar{x}_{i_s}$ - note that $i_s$ will be the value of the index in the input of the function \emph{{\bf Dbs.m}}- and the elements of the mesh differentiation matrix corresponding to $x_i$ are given in terms of those of the stencil matrix corresponding to $i_s$ by $D_s (i, j=[j_1, j_2,...j_M])=\bar{D}_s (i_s, m=[1,...M])$.

For the particular case in which $M=N$, the stencil  and the mesh coincide and, therefore, the entire mesh differentiation matrix is the same as the stencil matrix ($D_s=\bar{D}_s$). An example of this situation is provided by the well known Chebyshev differentiation  matrices [Trefethen (2001)] which are constructed by Lagrange interpolation on a mesh of Gauss-Lobatto points.  In the general case $M < N$,  it is natural to take $M$ odd, so that the stencil has a central point $\bar{x}_{M_c}$ with $M_c=(M+1)/2$. Thus, if $ M_c \le i \le N-M_c+1$,  $x_i$ is an interior mesh point and we take centered at it its associated stencil, i.e.,  $i_s=M_c$ and $\bar{x}_{M_c}=x_i$.  However, mesh points near the boundaries, defined  as  those $x_i$ with $\quad 
 i < M_c$ (left boundary) or $\quad  i > N-M_c+1$ (right boundary), have equal associated stencils formed by $M$ mesh points contiguous to each boundary. More explicitly, the following three cases are distinguished:
\begin{enumerate}
\item $\bf{1\le i < M_c}$, i.e., point $x_i$ is near the left boundary. Its associated stencil is defined to have  $i_s=i$  and elements
\begin{equation}
 \left(\bar{x}_{1}, \bar{x}_{2},...\bar{x}_{M}\right)=\left(x_{1}, x_{2},...x_{M}\right)\quad (\bar{x}_{i_s=i}=x_i ).  \label{left}
\end{equation}  
The differentation matrix elements corresponding  to $x_i$ are given by
\begin{equation}
D_s (i, j=[1,..M])= \bar{D}_s (i_s, m=[1,...M]) \quad (i \le M_c) . \label{Dsleft}
\end{equation}  

\item $\bf{M_c \le i \le N-M_c}$, i.e., point $x_i$ is a mesh interior point.  Its associated stencil is defined to have $i_s=M_c$ and elements
\begin{equation}
 \left(\bar{x}_{1}, \bar{x}_{2},...\bar{x}_{M}\right)=\left(x_{i-M_c+1}, x_{i-M_c+2},...x_{i+M_c-1}\right)\quad (\bar{x}_{i_s=M_c}=x_i), \label{interior}
\end{equation}  
The differentation matrix elements corresponding  to $x_i$ are given by
\begin{equation}
D_s (i, j=[i-Mc+1,...i+M_c-1])= \bar{D}_s (i_s,m=[1,...M]) \quad [M_c \le  i \le N-M_c] \label{Dsinterior}
\end{equation}

\item $\bf{ N-M_c < i \le N}$, i.e., the point $x_i$ is near the right boundary. Its associated stencil is defined to have $i_s=i+M-N$ and elements
\begin{equation}
\left(\bar{x}_{1}, \bar{x}_{2},...\bar{x}_{M}\right)=\left(x_{N-M+1}, x_{N-M+2},...x_{N}\right)\quad (\bar{x}_{i_s=i+M-N}=x_i). \label{right}
\end{equation}
 The differentation matrix elements corresponding  to $x_i$ are given by  
\begin{equation}
D_s (i, j=[N-M+1,..N])= \bar{D}_s (i_s, m=[1,...M]) \quad ( N- M_c+1 < i  \le N) . \label{Dsleft}
\end{equation}  
\end{enumerate}

Next we list the MATLAB function \emph{{\bf Dsmesh.m}} which, driving the function \emph{{\bf Dbs.m}} introduced in the previous section, computes the differentiation matrices   $D_s$  $(s=1,2,3)$ using stencils of any given length $M$ in an arbitray mesh of $N$ points stored in the vector ${\bf x}=(x_1,...x_N)$. 

\pagebreak

\lstset{language=Matlab, breaklines=true, basicstyle=\footnotesize}
\lstset{numbers=left, numberstyle=\tiny, stepnumber=1, numbersep=-2pt}
\begin{lstlisting}[frame=single]
 function [D_1 D_2 D_3]=Dsmesh(N, x, M)
  D_1=sparse(N,N); D_2=sparse(N,N); D_3=sparse(N,N);
  for i=1:N,      
    if M==N,
        xb=x; i_s=i;
        [D_1(i,1:M) D_2(i,1:M) D_3(i,1:M)]=Dbs(i_s,xb,M);
    else    
    M_c=(M+1)/2;
    if i < M_c,
        xb=x(1:M); i_s=i;
        [D_1(i,1:M) D_2(i,1:M) D_3(i,1:M)]=Dbs(i_s,xb,M);
    end 
    if i >= M_c && i<= N-M_c,
        xb=x((i-M_c+1):(i+M_c-1)); i_s=M_c
        [D_1(i,(i-M_c+1):(i+M_c-1)) D_2(i,(i-M_c+1):(i+M_c-1))...
            D_3(i,(i-M_c+1):(i+M_c-1))]=Dbs(i_s,xb,M);
    end       
    if i > N-Mc,
        xb=x((N-M+1):N); i_s=i+M-N;
        [D_1(i,(N-M+1):N) D_2(i,(N-M+1):N) D_3(i,(N-M+1):N)]=Dbs(i_s,xb,M);
    end    
   end
  end
\end{lstlisting}

 \section*{References}

[1] J.P. BERRUT and L. N. TREFETHEN, \emph{ Barycentric Lagrange Interpolation}, SIAM Review, 46 (2004), pp. 501-517.

[2] B. FORNBERG,\emph{A practical guide to pseudospectral methods}, Cambridge University Press, 1996.

[3] P. HENRICI, \emph{Essentials of numerical analysis}, Wiley, New York, 1982.

[4] R. J. LEVEQUE, \emph{Finite differences methods for ordinary and partial differential equations}, Society for Industrial and Applied Mathematics, Philadelphia, 2007.

[5] H. RUTIHAUSER, \emph{Vorlesungen uber numerische Mathematik, Vol. 1}, Birkhauser, Basel, Stuttgart, 1976; English translation, \emph{Lectures on Numerical Mathematics}, Walter Gautschi, ed., Birkhauser, Boston, 1990.

[6]  L. N. TREFETHEN, \emph{ Spectral methods with Matlab}, Society for Industrial and Applied Mathematics, Philadelphia, 2001.

[7]  J.A.C. WEIDEMAN and S.C. REDDY,  \emph{ A MATLAB differentiation matrix suite}, ACM Transactions on Mathematical Software, Vol. 26, No. 4, December 2000, Pages 465-519.

[8] B.D. WELFERT, \emph{ Generation of pseudospectral differentiation matrices I}, SIAM J. Numer. Anal., 34 (1997), pp. 1640-1657.

\end{document}